\documentclass[11pt]{amsart}
\usepackage{amsmath,amssymb,latexsym,amsthm,enumerate}

\theoremstyle{plain}

\theoremstyle{definition}

\theoremstyle{remark}

\newcounter{numpar}[section]
\newcommand*{\newpar}[1]{\refstepcounter{numpar}\medskip\textbf{\thenumpar.#1}}

\newcommand*{\dd}{\mathrm d}     

\newcommand*{\cF}{\mathcal F}

\newcommand*{\bbR}{\mathbb R}

\newcommand*{\EE}{\mathsf E}
\newcommand*{\PP}{\mathsf P}
\newcommand*{\QQ}{\mathsf Q}

\DeclareMathOperator{\sgn}{sgn}

\begin{document}
\title{A note on a paper by Wong and Heyde}

\author{Aleksandar Mijatovi\'{c}}
\address{Department of Statistics, University of Warwick, UK}
\email{a.mijatovic@warwick.ac.uk}

\author{Mikhail Urusov}
\address{Institute of Mathematical Finance, Ulm University, Germany}
\email{mikhail.urusov@uni-ulm.de}


\keywords{Local martingales vs. true martingales; stochastic exponentials}

\subjclass[2000]{60G44, 60G48, 60H10, 60J60}


\begin{abstract}
In this note we re-examine the analysis of the paper ``On the martingale property of stochastic exponentials''
by B.~Wong and C.C.~Heyde~\cite{WongHeyde:04}. 
Some counterexamples are presented
and alternative formulations are discussed. 
\end{abstract}

\maketitle

\section{Introduction}
\label{sec:intro} In~\cite{WongHeyde:04} the authors announce very
general results about the martingale property of exponential local
martingales. Since the subject matter of the paper is important, 
it is necessary to examine it critically.
In Section~\ref{sec:setting} we describe the setting of~\cite{WongHeyde:04}.
In Section~\ref{sec:discussion} we re-examine the analysis
in~\cite{WongHeyde:04}
and discuss an alternative formulation of their Corollary~2.
Section~\ref{sec:counter} contains counterexamples. 

\section{The setting in~\cite{WongHeyde:04}}
\label{sec:setting}
The main object of study in~\cite{WongHeyde:04} is
the stochastic exponential
\begin{eqnarray}
\label{eq:setting1}
Z_X(t) & = & \exp\left\{\int_0^tX(u)\cdot\dd W(u)-\frac{1}{2}\int_0^t\|X(u)\|^2\,\dd
u\right\},\quad t\in[0,\infty),\\
Z_X(0) & = & 1,
\nonumber
\end{eqnarray}
where $X(t)\in\bbR^d$ is a $d$-dimensional
$(\cF_t)$-progressively measurable process and
$W$
is a $d$-dimensional Brownian motion
defined on a given probability space
$(\Omega,\cF,\PP)$.
The filtration
$(\cF_t)_{t\in[0,\infty)}$
is assumed to be generated by
$W$
and augmented to satisfy the usual conditions
(see paragraph two on page 656 in~\cite{WongHeyde:04}
for these assumptions).

The stopping time
$\tau^{M_X}$
is defined by
\begin{equation}
\label{eq:setting2}
\tau^{M_X}=\lim_{N\to\infty}\tau^{M_X}_N,\quad\text{where}\quad
\tau^{M_X}_N=\inf\left\{t\in[0,\infty)\colon\int_0^t\|X(u)\|^2\,\dd u \geq N\right\},
\end{equation}
with the usual convention that $\inf\emptyset=\infty$
(see equation~(2)~on page~656 in~\cite{WongHeyde:04}).
The non-decreasing adapted process
$M_X$,
defined by
$$
M_X(t) := \int_0^t\|X(u)\|^2\,\dd u,\quad\text{for}\quad t\in[0,\infty),
$$
is left-continuous by the monotone convergence theorem.
On the event
$\left\{\tau^{M_X}<\infty\right\}$,
at time
$\tau^{M_X}$
the process
$M_X$
either jumps to infinity
if
$M_X\left(\tau^{M_X}\right)<\infty$
or tends to infinity continuously
(i.e. assumes arbitrarily large values just before
$\tau^{M_X}$)
if
$M_X\left(\tau^{M_X}\right)=\infty$.
A precise definition of the stochastic exponential
$Z_X$
in~\eqref{eq:setting1}
can now be given by
\begin{eqnarray}
\label{eq:setting3}
Z_X(t)  & :=  &  \exp\left\{\int_0^{t\wedge\tau^{M_X}}X(u)\cdot\dd W(u)-\frac{1}{2}\int_0^{t\wedge\tau^{M_X}}\|X(u)\|^2\,\dd u\right\},
\quad\text{for}\quad t\in[0,\infty),
\end{eqnarray}
where we set
$Z_X(t) = 0$
on
$\left\{\tau^{M_X}\leq t,\>M_X\left(\tau^{M_X}\right)=\infty\right\}$.
Note that
the stochastic integral
in~\eqref{eq:setting3}
is well-defined
$\PP$-a.s. on
$\left\{\tau^{M_X}<\infty,\, M_X\left(\tau^{M_X}\right)<\infty \right\}$
for every
$t\in[0,\infty)$,
and hence, the process
$Z_X$
is continuous and
takes strictly positive values on the event
$\left\{\tau^{M_X}=\infty\right\}\cup\left\{\tau^{M_X}<\infty,\, M_X\left(\tau^{M_X}\right)<\infty \right\}$.
Furthermore the stochastic exponential
$Z_X$
is a continuous local martingale (the continuity
at
$\tau^{M_X}$
on the event
$\{\tau^{M_X}<\infty,\,M_X\left(\tau^{M_X}\right)=\infty\}$
follows by the
Dambis-Dubins-Schwarz theorem, see e.g.~\cite[Ch.~V, Th.~1.6]{RevuzYor:99}).

On page 656 of~\cite{WongHeyde:04},
in the line following the formula containing the definition of
$\tau_N^{M_X}$, the authors stipulate that their process $X$ is stopped at $\tau^{M_X}$,
which is rather confusing because it implies $\PP\left(\tau^{M_X}<\infty,M_X\left(\tau^{M_X}\right)<\infty\right)=0$
and thus restricts the generality without being essential for the setting and what follows
(see, however, a more precise and detailed discussion on a related point in item~\ref{item:disc5} of Section~\ref{sec:discussion}).
Note that the event $\{\tau^{M_X}<\infty,M_X(\tau^{M_X})<\infty\}$ can in general be the entire space $\Omega$
(e.g. take $X$ deterministic).

Finally, it should be noted that the authors of~\cite{WongHeyde:04}
work on a finite time interval $[0,T]$
(see e.g. paragraph two on page~656 in~\cite{WongHeyde:04}),
while the setting introduced above is the infinite time horizon setting.
This difference does not play a role for the exposition below
but lets us quote many formulas from~\cite{WongHeyde:04} 
exactly as they are stated there without introducing inconsistency. 
(There are notational inconsistencies in~\cite{WongHeyde:04}
related to this point: e.g. they define $\tau_N^{M_X}$ exactly as in~\eqref{eq:setting2},
but if one works on a finite time interval $[0,T]$ and the process $X(t)$ is given for $t\in[0,T]$,
one should have used either the definition
$$
\inf\left\{t\in[0,T]\colon\int_0^t\|X(u)\|^2\,\dd u \geq N\right\}
$$
or the definition
$$
\inf\left\{t\in[0,T]\colon\int_0^t\|X(u)\|^2\,\dd u \geq N\right\}\wedge T
$$
for $\tau_N^{M_X}$.)

\section{Discussion of Section~3 in~\cite{WongHeyde:04}}
\label{sec:discussion}

In this section
we reinspect Proposition~1 in~\cite{WongHeyde:04}, pointing to problems in its formulation and
proof. This has consequences for the rest of the paper~\cite{WongHeyde:04}. The formulation and
the proof of the main result, Theorem~1 in~\cite{WongHeyde:04}, both rely on
Proposition~1.
Further, in Definition~1 in~\cite{WongHeyde:04} the
central concept of a ``candidate measure'' is introduced. It is
implicitly assumed throughout~\cite{WongHeyde:04} that the defined
object exists and is unique.
However the ``candidate measure'' may not in fact exist, and
if it does, it may not be unique.
Thus Corollaries~1 and~2 in~\cite{WongHeyde:04}, which in their formulation use the notion of
the ``candidate measure'', are invalid.
Corollary~2 in~\cite{WongHeyde:04} admits an obvious well-posed reformulation,
but the resulting statement, given in Corollary~2' below, is also invalid 
(see Section~\ref{sec:counter} for a counterexample).

\newpar{}
\label{item:disc1}
We first discuss
Proposition~1
in~\cite{WongHeyde:04},
which plays the key role in~\cite{WongHeyde:04}
as it is used in the formulation and
applied in the proof of the main result of~\cite{WongHeyde:04}
(Theorem~1).
We start by quoting Proposition~1 of~\cite{WongHeyde:04} (see Section~3, page~657).

\smallskip
\noindent \textbf{Proposition 1 of [WH04].} \textit{Consider a
$d$-dimensional
$\cF_t$-progressively measurable process
$X(t)=\xi(W(\cdot),t)$
defined possibly up to the explosion time
$\tau^{M_X}$
defined by~\eqref{eq:setting2}.
Then there will also exist a
$d$-dimensional
$\cF_t$-progressively measurable process
$Y(t)=\xi(W(\cdot)+\int_0^\cdot Y(u)\,\dd u,t)$,
defined possibly up to the explosion time
$\tau^{M_Y}$,
with
$$
\tau^{M_Y}=\lim_{N\to\infty}\tau^{M_Y}_N,
$$
where
\begin{eqnarray*}
M_Y(t) & = & \int_0^t\|Y(u)\|^2\,\dd u,\\
\tau^{M_Y}_N & = & \inf\left(t\in[0,\infty)\colon\int_0^t\|Y(u)\|^2\,\dd u \geq N\right).
\end{eqnarray*}
}

The formulation of Proposition 1 is misleading.
Firstly, the formula ``$Y(t)=\xi(W(\cdot)+\int_0^. Y(u)\,\dd u,t)$'' cannot be a definition of $Y$.
This is an equation in
$Y$.
Secondly,
the statement in
the proposition that
``... the process
$Y(t)=\xi(W(\cdot)+\int_0^. Y(u)\,\dd u,t)$
defined possibly
up to the explosion time $\tau^{M_Y}$ ...''
is followed by the definition
of $\tau^{M_Y}$, which is given in terms of $Y$
that has not yet been defined.

Let us now analyse the proof of Proposition~1 in~\cite{WongHeyde:04}
in the hope that it will shed light on its formulation.
The proof operates with a process $X_N$,
which is not introduced in~\cite{WongHeyde:04}.
However, the formula ``$Z_X(t\wedge\tau^{M_X})=Z_{X_N}(t)$''
in the first line of the proof
(which should read as ``$Z_X(t\wedge\tau_N^{M_X})=Z_{X_N}(t)$'',
as supported by what follows and because the authors of~\cite{WongHeyde:04}
refer to their Lemma~1 in the second line of the proof)
makes it evident that they mean $X_N(t)=X(t)I(t\le\tau_N^{M_X})$.
The authors of~\cite{WongHeyde:04} define a measure $\QQ_N$,
equivalent to $\PP$, by $\QQ_N(A) = \EE_\PP[Z_{X_N}(T)I(A)]$ for all $A\in\cF_T$.
(Here we corrected another misprint: in~\cite{WongHeyde:04} they write
``$\QQ_N\left(X_N\in A\right) = \EE_\PP\left[Z_{X_N}(T)I(X_N\in A)\right]$ for all $A\in\cF_T$''.)
Further they define a $d$-dimensional $\QQ_N$-Brownian motion $W^{\QQ_N}$
by the formula
$$
W^{\QQ_N}(t)=W(t)-\int_0^t X_N(u)\,\dd u.
$$
Then, in line~2 on page~658 in~\cite{WongHeyde:04} the identity
\begin{eqnarray}
\label{eq:discussion1}
X_N(t) & = & \xi\left(W^{\QQ_N}(\cdot)+\int_0^.X_N(u)\,\dd u,t\right)
\text{ on }\{t\le\tau_N^{M_X}\}
\end{eqnarray}
is stated, which is correct. However, all that follows in the proof of Proposition~1
has problems.
It is stated in~\cite{WongHeyde:04} that~\eqref{eq:discussion1} ``demonstrates the existence up to $\tau_N^{M_Y}$ of''
\begin{eqnarray}
\label{eq:discussion2}
Y_N(t) & = & \xi\left(W(\cdot)+\int_0^.Y_N(u)\,\dd u,t\right)
\end{eqnarray}
(see page~658, line~4 in~\cite{WongHeyde:04}).
Firstly, it is not clear how to understand the words ``up to $\tau_N^{M_Y}$''
since
$\tau_N^{M_Y}$ is defined through $Y$ in the formulation of Proposition~1,
while $Y$ is still undefined.
Secondly, this statement is incorrect because~\eqref{eq:discussion1}  is just an identity
that holds for the particular processes $X_N$ and $W^{\QQ_N}$,
while~\eqref{eq:discussion2}  is  an equation in $Y_N$, where $W$ is the given
initial Brownian motion under $\PP$.
Using an argument, similar to this one in~\cite{WongHeyde:04}, one can conclude
that the existence of a weak solution of a stochastic differential equation
``demonstrates the existence of'' a strong
solution of the same equation, which, however, is false as is well-known;
see, e.g., \cite[Ch.~5, Ex.~3.5]{KaratzasShreve:91} or item~\ref{item:disc5.5} of this section.
Furthermore, even if this transition from~\eqref{eq:discussion1} to~\eqref{eq:discussion2}
were in order,
one would not be able to take limits as $N\to\infty$
as suggested in~\cite{WongHeyde:04} (page~658, line~6)
because nothing is said about the uniqueness of $Y_N$
satisfying~\eqref{eq:discussion2} nor about the consistency properties
of the ``solutions'' $Y_N$ of~\eqref{eq:discussion2}  for different $N$.
One must conclude that the proof of Proposition~1 in~\cite{WongHeyde:04} is invalid,
in whichever way  
one interprets the statement.\\

This in turn invalidates the main result: Theorem~1 on page~658 in~\cite{WongHeyde:04}
is misleading since its formulation and proof use the process $Y$
from Proposition~1 of~\cite{WongHeyde:04}.

\newpar{}
\label{item:disc2}
In Definition~1
on page 660 of~\cite{WongHeyde:04},
which we now quote,
the authors
``define''
the measure
$\QQ^C$
as follows.

\smallskip
\noindent \textbf{Definition 1 of [WH04].}
\textit{A \emph{candidate measure}
$\QQ^C$,
corresponding to the process
$X(t)$ defined in Proposition~1 on the measure
$\PP$, is a measure corresponding to which
\begin{equation}
\label{eq:discussion6}
X(t)=\xi\left(W^{\QQ^C}(\cdot)+\int_0^\cdot X(u)\,\dd u,t\right)
\end{equation}
is defined (possibly up to the explosion time $\tau^{M_X}$),
with $W^{\QQ^C}$ a $\QQ^C$-Brownian motion.}

\smallskip
This ``definition'' is unclear, regarding both existence (of $\QQ^C$ and $W^{\QQ^C}$)
and uniqueness.  In~\cite{WongHeyde:04} the authors say that $\QQ^C$ is well-defined by the
analysis in Proposition~1, but this argument is invalid as discussed above.  Indeed, counterexamples in
items~\ref{item:counter1} and~\ref{item:counter2}
of Section~\ref{sec:counter} below show that both existence and uniqueness of $\QQ^C$
``defined'' in this way may fail.  This in turn invalidates the next result, Corollary~1, on
page~660 in~\cite{WongHeyde:04}.

\newpar{}
\label{item:disc3}
We now turn our attention to
Corollary~2 on page~661 in~\cite{WongHeyde:04},
which the authors formulate as follows.

\smallskip
\noindent \textbf{Corollary 2 of [WH04].} \textit{Assume that
$X(t)$
is the unique weak solution up to the explosion time
$\eta_X$
of the functional SDE}
\begin{equation}
\label{eq:discussion3}
\dd X(t) = \mu(X,t)\,\dd t+\sigma(X,t)\cdot \dd W(t)
\end{equation}
\textit{with initial value
$X(0)$
and
$\mu(x,t)\in\bbR^d,$
$\sigma(x,t)\in\bbR^{d\times r},$
with
$\mu(x,t), \sigma(x,t)$
progressively measurable functionals.
Then
\begin{eqnarray}
\label{eq:discussion4}
\EE_\PP[Z_X(T)] = \QQ^C(\eta^X>T),
\end{eqnarray}
where
\begin{equation}
\label{eq:discussion5}
\dd X(t) = (\mu(X,t)+\sigma(X,t)\cdot X(t))\,\dd t +\sigma(X,t)\cdot \dd W^{\QQ^C}(t).
\end{equation}
}

The first minor point here is that $X$ is $d$-dimensional and $W$ is $r$-dimensional,
while it is important in the definition of the process $Z_X$ that $X$ and $W$ have the same dimension.

The explosion time
$\eta^X$,
which appears in Corollary~2,
is defined in the last paragraph on page 660 in~\cite{WongHeyde:04}
by
\begin{eqnarray*}
\eta^X=\lim_{N\to\infty}\eta^X_N,\quad
\text{where}\quad
\eta^X_N=\inf\Big\{t\in[0,\infty)\colon\sup_{i=1,\ldots,d}|X_i(t)|\geq N\Big\},
\end{eqnarray*}
$X_i(t)$, $i=1,\ldots,d$ are components of $X(t)$.
Let us add at this point that both $\tau^{M_X}$ and $\eta^X$
are termed ``the explosion time'' in~\cite{WongHeyde:04}
(see e.g. the above formulations of Proposition~1 and Definition~1
quoted from~\cite{WongHeyde:04}),
which is confusing because these stopping times can be different even in the setting of Corollary~2
(e.g. take an appropriate deterministic~$X$).
In our paper only $\eta^X$ is called ``the explosion time''
with the exception of the statements that we quote from~\cite{WongHeyde:04}.\\
\indent We conclude that Corollary 2 of [WH04] is also invalid, as it is unclear what the
measure $\QQ^C$ represents. \\

\newpar{}
\label{item:disc4}
We now seek a well-posed reformulation of Corollary 2 of [WH04].  There is a natural candidate,
as follows (though this is still incorrect, as we show next). 
Since $W^{\QQ^C}$ is assumed to be a $\QQ^C$-Brownian motion
and SDE~\eqref{eq:discussion5} is announced to hold,
it is natural to suggest the following: 

\smallskip
\noindent \textbf{Corollary 2'.} \textit{Let
$X$
be a unique in law possibly explosive weak solution of the SDE
\begin{eqnarray}
\label{eq:discussion7}
\dd X(t) = \mu(X,t)\,\dd t+\sigma(X,t)\cdot \dd W(t)
\end{eqnarray}
on some filtered probability space $(\Omega,\mathcal F,(\mathcal F_t),\PP)$
with initial value
$X(0)$
and
$\mu(x,t)\in\bbR^d,$
$\sigma(x,t)\in\bbR^{d\times d},$
where
$\mu(x,t), \sigma(x,t)$
are progressively measurable functionals.
Consider the process
\begin{eqnarray}
\label{eq:discussion8}
Z_X(t)  & =  &  \exp\left\{\int_0^{t\wedge\eta^{X}}X(u)\cdot\dd W(u)-\frac{1}{2}\int_0^{t\wedge\eta^{X}}\|X(u)\|^2\,\dd u\right\},
\quad\text{for}\quad t\in[0,\infty),
\end{eqnarray}
where we set
$Z_X(t)=0$
for
$t\geq\eta^X$
on the event
$\{\int_0^{\eta^X}\|X(u)\|^2\,\dd u=\infty\}$.
Assume further that
$\widetilde X$ is a unique in law possibly explosive weak solution of the SDE
\begin{eqnarray}
\label{eq:discussion9}
\dd \widetilde X(t) = (\mu(\widetilde X,t)+\sigma(\widetilde X,t)\cdot \widetilde X(t))\,\dd t +\sigma(\widetilde X,t)\cdot \dd \widetilde W(t)
\end{eqnarray}
on some filtered probability space $(\widetilde\Omega,\widetilde{\mathcal F},(\widetilde{\mathcal F}_t),\widetilde\PP)$
with the same initial value $X(0)$.
Then
$$
\EE_\PP[ Z_X(T)] = \widetilde \PP(\eta^{\widetilde X}>T),
$$
where
$\eta^{\widetilde X}$
is the explosion time of
$\widetilde X$.}

\smallskip
Let us point out the difference between Corollary~2 and Corollary~2'
in that in the latter the existence of unique in law
weak solution is assumed for each SDE separately,
possibly on different probability spaces,
while in the former both measures $\PP$ and $\QQ^C$
are stated to be on the same space and the process
$X$ is claimed to solve the two SDEs under the two measures respectively.

\newpar{}
\label{item:disc5}
Before we proceed with Corollary~2' let us point out a further inconsistency
in the formulation of Corollary~2 in~\cite{WongHeyde:04},
which is also fixed in the formulation of Corollary~2'.
In the formulation of Corollary~2
the authors of~\cite{WongHeyde:04} go beyond their setting.
Namely, allowing
$X$
to be explosive is inconsistent with their
definition of $\tau^{M_X}$ and
stipulation that the process $X$
is stopped at $\tau^{M_X}$.
This is particularly relevant if e.g. we have
\begin{eqnarray}
\label{eq:discussion10}
\eta^X<\infty \quad \PP\text{-a.s.}
&\text{ and }&
\int_0^{\eta^X}\|X(u)\|^2\,\dd u<\infty\quad \PP\text{-a.s.},
\end{eqnarray}
since in this case~\eqref{eq:setting2}
does not define
$\tau^{M_X}$
(unless it is specified what $X$ is after $\eta^X$;
note that it may happen that the limit $\lim_{t\uparrow\eta^X}X(t)$ does not exist).
It is easy to see that~\eqref{eq:discussion10} is indeed possible
(take e.g. an appropriate deterministic~$X$,
which corresponds to zero matrix $\sigma$ in~\eqref{eq:discussion3};
in item~\ref{item:counter4} of Section~\ref{sec:counter}
we also give a stochastic example, where \eqref{eq:discussion10}~holds).
Thus, if \eqref{eq:discussion10}~holds
and the behaviour of $X$ after $\eta^X$ is not specified
(which is the case in the setting of Corollary~2 in~\cite{WongHeyde:04}),
then $\tau^{M_X}$ is undefined, $Z_X$ is also undefined
(the authors of~\cite{WongHeyde:04} define $Z_X$ via~\eqref{eq:setting3}; see Section~2 in~\cite{WongHeyde:04}),
and hence, the left-hand side of~\eqref{eq:discussion4} is undefined as well.

In the proof of Corollary~2 on page~661 in~\cite{WongHeyde:04}
the authors claim that $\eta^X=\tau^{M_X}$~a.s.\footnote{In the formula in lines~2 and~3
of the proof of Corollary~2 in~\cite{WongHeyde:04}
it is stated that~$\eta^X$ and $\tau^{M_X}$ have the same law
(note that the mentioned formula is claimed to hold for any $T$).
However, since it is clear that $\eta^X\leq \tau^{M_X}$~a.s.
(provided it is specified what $X$ is after $\eta^X$,
to be able to speak about~$\tau^{M_X}$),
this claim amounts to $\eta^X=\tau^{M_X}$~a.s.},
which is incorrect as we have just seen that $\tau^{M_X}$
may be undefined in the setting of Corollary~2 in~\cite{WongHeyde:04}.
Moreover, even if the behaviour of $X$ after $\eta^X$ were specified
(so that $\tau^{M_X}$ were well-defined),
then the claim $\eta^X=\tau^{M_X}$~a.s. would be also incorrect
(e.g. if \eqref{eq:discussion10}~holds and we specify
$X(t)=0$ for $t\ge\eta^X$, then $\tau^{M_X}=\infty>\eta^X$~$\PP$-a.s.).
In order to define $Z_X$ in the setting of Corollary~2 in~\cite{WongHeyde:04}
one needs to use formula~\eqref{eq:discussion8},
in which case no problems arise and the behaviour of $X$
after $\eta^X$ in not essential at all.

\newpar{}
\label{item:disc5.5}
Furthermore, it should be emphasised that
Corollary~2, as stated in~\cite{WongHeyde:04},
goes beyond the setting of~\cite{WongHeyde:04}
also in another respect, and hence,
even if there were no issues with Proposition~1, Definition~1
and other issues with Corollary~2 discussed above,
Corollary~2 in~\cite{WongHeyde:04}
could not be proved as claimed.
Recall that a \emph{solution} (or a \emph{weak solution}) of the SDE
$$
\dd X(t)=\mu(X,t)\,\dd t+\sigma(X,t)\,\dd W(t),\quad X(0)=x_0,
$$
is a pair of adapted processes $(X,W)$
on some filtered probability space
$(\Omega,\cF,(\cF_t),\PP)$ such that $W$
is an $(\cF_t)$-Brownian motion,
$$
\int_0^t(|\mu(X,u)|+\sigma^2(X,u))\,\dd u<\infty\quad\PP\text{-a.s.},\quad t\in[0,\infty),
$$
and
$$
X(t)=x_0+\int_0^t\mu(X,u)\,\dd u+\int_0^t\sigma(X,u)\,\dd W(u)\quad\PP\text{-a.s.},\quad t\in[0,\infty).
$$
(For notational simplicity, we consider only one-dimensional $X$ and $W$ and define only a non-explosive solution here
because this is all that we need for the argument with Tanaka's SDE below.)
Note that the filtered probability space
$(\Omega,\cF,(\cF_t),\PP)$
in this definition may differ from the one described
in Section~\ref{sec:setting}.
A \emph{strong solution}
$(X,W)$
of this SDE
is a solution such that the
process $X$
is adapted to the filtration generated by the Brownian motion $W$
(see~\cite[Ch.~IX, \S~1]{RevuzYor:99} for the employed terminology).
It is easy to show that Tanaka's SDE
\begin{equation}
\label{eq:discussion11}
\dd X(t) = \sgn X(t)\,\dd W(t),\quad X_0=0,
\end{equation}
where
$$
\sgn x=\begin{cases}
1&\text{if }x>0,\\
-1&\text{if }x\le0,
\end{cases}
$$
has a unique in law weak solution
(and, moreover, $X$ is a Brownian motion by L\'evy's characterisation theorem
for any solution $(X,W)$ of~\eqref{eq:discussion11}).
However, there exists no strong solution of~\eqref{eq:discussion11}:
for any solution
$(X,W)$
we have
$$
W(t)=\int_0^t\sgn X(u)\,\dd X(u)\quad\PP\text{-a.s.}\quad\text{for all}\quad t\geq0,
$$
and hence,
by~\cite[Ch.~VI, Cor.~2.2]{RevuzYor:99},
the filtration generated by $W$ coincides with that generated by $|X|$,
which is strictly smaller than the filtration generated by $X$, as $X$ is a Brownian motion.
Thus, $X$
cannot be adapted to the filtration generated by
$W$.
This argument implies that a solution of Tanaka's SDE
cannot be expressed as
$X(t)=\xi(W(\cdot),t)$,
for a progressively measurable functional $\xi$,
and hence
does not satisfy the assumptions of Proposition~1 and Theorem~1 in~\cite{WongHeyde:04}.
To summarize the last point, even if all other results in~\cite{WongHeyde:04} were
beyond reproach, the weak existence, which is assumed in Corollary~2 in~\cite{WongHeyde:04},
would be an insufficient assumption to support the conclusions of Corollary 2 by using their method.
The authors of~\cite{WongHeyde:04} should have assumed existence of a strong solution
of~\eqref{eq:discussion3} in Corollary~2.

\newpar{}
\label{item:disc6}
The discussion above leads to the question of whether Corollary 2' holds.
Or, at least, whether such a statement holds under the stronger assumptions
that $X$ and $\widetilde X$ are pathwise unique strong solutions of
SDEs \eqref{eq:discussion7} and~\eqref{eq:discussion9}.
A counterexample in item~\ref{item:counter3} of Section~\ref{sec:counter} shows that the answer is negative.
Moreover, $X$ and $\widetilde X$ are pathwise unique strong solutions of those SDEs in that counterexample.

%

\section{Counterexamples}
\label{sec:counter}
\newpar{}
\label{item:counter1}
We start with two counterexamples to Definition~1 in~\cite{WongHeyde:04}.
Let us take $d=1$, fix a finite time horizon $T\in(0,\infty)$,
and consider $\Omega=C([0,T],\bbR)$ the space of continuous functions $[0,T]\to\bbR$.
Let $W$ be the coordinate process on~$\Omega$,
$\PP$~the Wiener measure,
$(\cF_t)_{t\in[0,T]}$ the filtration generated by $W$
and augmented to satisfy the usual conditions, and $\cF=\cF_T$,
so that we are in the setting of~\cite{WongHeyde:04}.

First we show that the measure $\QQ^C$ in Definition~1 in~\cite{WongHeyde:04} may not be unique.
Indeed, take $\xi(\cdot,\cdot)\equiv0$, so that $X\equiv0$ as well
(recall that $X$ is defined by the formula $X(t)=\xi(W(\cdot),t)$).
Note that, for any $\lambda\in\bbR$, the process
$$
W^\lambda(t)=W(t)-\lambda t,\quad t\in[0,T],
$$
is an $(\cF_t,\PP^\lambda)$-Brownian motion, where the measure $\PP^\lambda$ is given by
$$
\frac{d\PP^\lambda}{d\PP}=\exp\left\{\lambda W(T)-\frac{\lambda^2}2 T\right\}.
$$
Clearly, any measure $\PP^\lambda$ (and, in fact, many other measures)
can be considered as $\QQ^C$ because \eqref{eq:discussion6} is satisfied with $W^{\QQ^C}=W^\lambda$,
which is a $\PP^\lambda$-Brownian motion as required in Definition~1.

\newpar{}
\label{item:counter2}
Now we show that the measure $\QQ^C$ in Definition~1 in~\cite{WongHeyde:04} may not exist.
We consider the filtered probability space as above and a strictly increasing continuous function
$f\colon\bbR\to\bbR$
with
\begin{equation}
\label{eq:counter1}
\lim_{x\to-\infty}f(x)=1\text{ and }\lim_{x\to\infty}f(x)=2.
\end{equation}
Let us define a progressively measurable functional $\xi(\omega,t)$, $\omega\in\Omega(=C([0,T],\bbR))$, $t\in[0,T]$,
by the formula
\begin{equation}
\label{eq:counter2}
\xi(\omega,t)=\begin{cases}
\frac{f(\omega(t))}{T-t}&\text{if }t\in[0,T),\\
0&\text{if }t=T,
\end{cases}
\end{equation}
which gives us the process $X(t)=\xi(W(\cdot),t)$, $t\in[0,T]$.

Let us prove that there exists no measure $\QQ^C$ satisfying Definition~1 in~\cite{WongHeyde:04}.
Since $X(t)=f(W(t))/(T-t)$, $t\in[0,T)$,
and $f$ is strictly increasing,
trajectories of $X$ determine trajectories of $W$ uniquely.
In particular, if \eqref{eq:discussion6} holds,
the process $W^{\QQ^C}$ should satisfy
$$
W^{\QQ^C}(t)=W(t)-\int_0^t X(u)\,\dd u,\quad t\in[0,T).
$$
It follows that
$$
\lim_{t\uparrow T}W^{\QQ^C}(t)=-\infty\text{\emph{ for any }}\omega\in\Omega
$$
(recall \eqref{eq:counter1},~\eqref{eq:counter2} and note that $W(t)\to W(T)\in\bbR$ as $t\uparrow T$ \emph{for any} $\omega\in\Omega$
because $W$ is the coordinate process on the space of continuous functions),
hence there does not exist a measure $\QQ^C$ on $(\Omega,\cF)$ such that $W^{\QQ^C}$ is a $\QQ^C$-Brownian motion.

\newpar{}
\label{item:counter3}
We proceed with a counterexample to Corollary~2'.
Let
$\mu(x)=|x|^\alpha$,
for any fixed
$\alpha>3$,
and
$\sigma(x)\equiv 1$.
The process
$X(t)$,
with the state space
$\bbR$
and starting value
$X(0)\in\bbR$,
can be defined as a strong solution of the SDE
\begin{equation}
\label{eq:counter3}
\dd X(t)=|X(t)|^\alpha\,\dd t+\dd\, W(t)
\end{equation}
up to the explosion time
$\eta^X$.
The existence of a strong solution up to
$\eta^X$ and pathwise uniqueness are
guaranteed by It\^o's existence and uniqueness theorem
since the coefficients of the SDE are locally Lipschitz
(see~\cite[Ch.~IX, Ex.~2.10]{RevuzYor:99}).
It follows from Example 3.1 of~\cite{MijatovicUrusov:10a}
that the process
\begin{eqnarray*}
Z_X(t)  & = & \exp\left\{\int_0^{t\wedge\eta^{X}}X(u)\,\dd W(u)
-\frac12\int_0^{t\wedge\eta^{X}}X^2(u)\,\dd u\right\},
\quad t\in[0,\infty)
\end{eqnarray*}
($Z_X(t)=0$ for $t\ge\eta^{X}$ on $\{\int_0^{\eta^X} X^2(u)\,\dd u=\infty\}$;
see, however, formula~\eqref{eq:counter5} below) is a martingale
(in fact, it is even a uniformly integrable martingale).
Hence, we have
$$
\EE_\PP[Z_X(T)] = 1 \quad\text{for all}\quad T\geq0.
$$
In this case, SDE~\eqref{eq:discussion9} has the form
\begin{eqnarray}
\label{eq:counter4}
\dd\widetilde X(t)=\left(|\widetilde X(t)|^\alpha+\widetilde X(t)\right)\,\dd t+\dd\widetilde W(t).
\end{eqnarray}
Its coefficients are locally Lipschitz and therefore there exists a pathwise unique strong solution up to the explosion time $\eta^{\widetilde X}$.
The process $\widetilde X$ explodes to $+\infty$ in finite time, which follows from Feller's
test for explosions (see \cite[Ch.~5, Th.~5.29 and Prop.~5.32]{KaratzasShreve:91}).
By time-homogeneity of SDE~\eqref{eq:counter4},
$$
\widetilde\PP\left(\eta^{\widetilde X}>T\right)<1\quad\text{for all}\quad T>0,
$$
which now contradicts the claim in Corollary~2'.

Note that since SDE~\eqref{eq:counter4}
has a pathwise unique strong solution, we can construct a solution
of this SDE on the same probability space that supports the
solution of SDE~\eqref{eq:counter3} with the same Brownian motion
$W$
as in~\eqref{eq:counter3}.
This means that the reason why Corollary 2' does not hold
is not due to the fact that the solutions of SDEs~\eqref{eq:discussion7}
and~\eqref{eq:discussion9} are allowed to exist on distinct probability spaces;
in fact this reason is deeper. See~\cite{MijatovicUrusov:10a} for more
details on this point.

\newpar{}
\label{item:counter4}
Finally, as promised in Section~\ref{sec:discussion},
we demonstrate that~\eqref{eq:discussion10} is possible.
Namely,~\eqref{eq:discussion10} holds in the example in item~\ref{item:counter3} of this section.
Indeed, by Feller's test for explosions, $\PP$-almost all
trajectories of
$X$
explode at
$+\infty$,
and hence,
$\eta^X<\infty$~$\PP$-a.s.
We now need to prove that
\begin{eqnarray}
\label{eq:counter5}
\int_0^{\eta^X}X_s^2\,\dd s<\infty\quad \PP\text{-a.s.}
\end{eqnarray}
in this example. The property~\eqref{eq:counter5}
is equivalent to
${Z_X(\eta^X)>0}$~$\PP$-a.s.,
which is in turn equivalent to the property
\begin{eqnarray}
\label{eq:counter6}
Z_X(\infty)>0\quad\PP\text{-a.s.}
\end{eqnarray}
(note that
$Z_X$
is stopped at
$\eta^X$). It remains to note that~\eqref{eq:counter6}
holds in the case when the process
$X$
is given by~\eqref{eq:counter3}
(with $\alpha>3$ as above),
which follows from Theorem~2.2 in~\cite{MijatovicUrusov:10b}.
Namely, condition~(II) in Theorem~2.2
in~\cite{MijatovicUrusov:10b} is satisfied.


\bibliographystyle{alpha}
\bibliography{refs}
\end{document}